\documentclass[12pt]{article}
\usepackage{amsmath}
\usepackage{amssymb}
\begin{document}

\newtheorem{theorem}{Theorem}[section]
\newtheorem{remark}[theorem]{Remark}
\newtheorem{mtheorem}[theorem]{Main Theorem}
\newtheorem{proposition}[theorem]{Proposition}
\newtheorem{lemma}[theorem]{Lemma}
\newtheorem{mlemma}[theorem]{Main-Lemma} 
\newtheorem{phlemma}[theorem]{Pigeon-Hole-Lemma}
\newtheorem{corollary}[theorem]{Corollary}
\newtheorem{example}[theorem]{Example}
\newtheorem{observation}[theorem]{Observation}
\newtheorem{construction}[theorem]{The $R$-module - Construction}
\newtheorem{definition}[theorem]{Definition}

\renewcommand{\labelenumi}{(\roman{enumi})}
\newcommand{\dach}[1]{\hat{\vphantom{#1}}}

\def\Z{{ \mathbb Z}}
\def\N{{ \mathbb N}}
\def\BB{ \mathbb B}
\def\BBB{ B_\BB}

\def\V{{\mathfrak V}}
\def\C{{\mathfrak C}}
\def\X{{\mathfrak X}}

\def\R{{\bf R}}
\def\D{\widehat{D}}
\def\A{\widehat{A}}
\def\Q{{\mathbb Q}}
\def\G{\hat{G}}
\def\T{{\cal T}}
\def\B{\widehat{B}}
\def\BC{\widehat{B_C}}
\def\BCC{\widehat{B_\C}}
\def\restr{\restriction}
\def\Aut{{\rm Aut\,}}
\def\Im{{\rm Im\,}}
\def\ker{{\rm ker\,}}
\def\inf{{\rm inf\,}}
\def\sup{{\rm sup\,}}
\def\Br{{\rm Br\,}}
\def\Yphi{Y_{[\phi]}}
\def\Ypsi{Y_{[\psi]}}
\def\Xphi{X_{\tilde{\phi}}}
\def\Xpsi{X_{\tilde{\psi}}}
\def\a{\alpha}
\def\abar{\overline{\alpha}}
\def\aa{{\bf a}}
\def\to{\rightarrow}
\def\arr{\longrightarrow}
\def\sigmaa{{\bf \Sigma_a}}

\def\End{{\rm End\,}}
\def\Ines{{\rm Ines\,}}
\def\Hom{{\rm Hom\,}}

\title{{\sc     Indecomposable almost free modules - the local case}
\footnotetext{This work is supported by the project 
No. G-0294-081.06/93 of the German-Israeli
Foundation for Scientific Research \& Development\\
AMS subject classification: 20K20, 20K26, 20K30, 13C10 \\
Key words and phrases: indecomposable modules of local rings,
$\aleph_1$-free modules of rank $\aleph_1$, realizing rings as endomorphism
rings\\
Sh 591 in Shelah's list of publications}
} 

\author{ R\"udiger G\"obel and Saharon Shelah}

\date{}

\maketitle

\begin{abstract}
Let $R$ be a countable, principal ideal domain which is
not a field and $A$ be a countable $R$-algebra which is free as an
$R$-module. Then we will construct an $\aleph_1$-free $R$-module $G$ of
rank $\aleph_1$ with endomorphism algebra End$_RG = A$. 
Clearly the result does not hold for fields.
Recall that
an $R$-module is $\aleph_1$-free if all its countable submodules are
free, a condition closely related to Pontryagin's theorem. This result
has many consequences, depending on the algebra $A$ in use. For
instance, if we choose $A = R$, then clearly $G$ is an indecomposable
`almost free' module. The existence of such modules was unknown for
rings with only finitely many primes like $R = \Z_{(p)}$, the
integers localized at some prime $p$. The result complements a
classical realization theorem of Corner's showing that any such
algebra is an endomorphism algebra of some torsion-free, reduced
$R$-module $G$ of countable rank. Its proof is based on new
combinatorial-algebraic techniques related with what we call {\it
rigid tree-elements} coming from a module generated over a forest of
trees. 
\end{abstract}

\section{\sc Introduction}

Let $R$ be a fixed countable, principal ideal domain which is not a
field.
An $R$-module $A$ is reduced for if  $\bigcap\nolimits_{s \in S}sA = 0$
where $S$ = $R \setminus \{0\}$ and $A$ is torsions-free if
$sa = 0$ ($s\in S, a\in A$) implies $ a = 0$. Note that $R$ is
reduced because $R$ is not a field.
We will consider  $R$-algebras $A$ which are torsion-free and reduced
as $R$-modules $A_R$. In particular this is the case when $A_R$ is
free.

Let $\kappa \leq\lambda$ be infinite cardinals. 
We are interested in $R$-modules of size $\lambda$ which are $\kappa$--free, 
which is the case when all its submodules of cardinality $< \kappa$ are free 
$R$-modules. 

{\em Can we find 
indecomposable $\kappa$--free $R$-modules of cardinality $\lambda$?}

We are mainly interested in the case when $\kappa = \lambda$ and
in particular when this cardinal is $\aleph_1$.

Such modules - by freeness - most likely want to decompose into non-trivial
direct sums and in fact, if $\lambda$ is a singular cardinal, then
by Shelah's \cite{S0} singular compactness theorem it follows that
such $R$-modules are free (hence very decomposable), 
this holds in particular for cardinals of 
cofinality $\omega$, e.g. for $\aleph_\omega$, a result due to Hill
\cite{Hi}, see Eklof, Mekler \cite {EM}.

On the other hand, the existence of non--free, $\aleph_1$--free
$R$-modules of cardinality $\aleph_1$ follows from Griffith \cite {Gr}, 
Hill \cite {Hi}, Eklof \cite {E1}, Mekler \cite {Me} and a result of 
Shelah's in Eklof \cite[p.82, Theorem 8.8]{E2}. By an induction it 
can be shown, that there are non--free $\aleph_n$--free modules of 
cardinality $\aleph_n$. A similar result for non--commutative groups 
is due to Higman \cite {Hg1,Hg2}.
The freeness-result at $\aleph_\omega$ illustrates that induction breaks 
down at $\aleph_{\omega + 1}$ and new techniques are needed to show that for
{\em certain cardinals $\lambda$ only} the existence of
non--free $\lambda$--free $R$-modules of cardinality $\lambda$ follows,
see Shelah, Magidor \cite{MS}.

\medskip
However only very little is known about algebraic properties of the 
non--free $\lambda$--free $R$-modules of cardinality
$\lambda$ and this is also the case when $\lambda = \kappa = \aleph_1$, 
see Eklof \cite {E1} and Eklof, Mekler \cite {EM}. The following problem 
is immediate.

\medskip
{\em Investigate the {\bf algebraic properties} of $\lambda$-free modules of 
cardinality $\lambda$}.
\medskip

The only earlier result known to us uses a construction from Shelah
\cite{S2} of non-separable groups \cite {E2,EM}
and is due to Eda  \cite {Ed}. He shows the existence of an
$\aleph_1$--free group  $G$  of cardinality $\aleph_1$ with
$\Hom(G,{\Z}) = 0$. In this paper we want to present new techniques which 
allow us to shed some more light on this problem.                                                                                                       

In order to work exclusively in ZFC we restrict ourself to 
$\kappa = \aleph_1$ and $\lambda \leq 2^{\aleph_0}$. Recall that
under negation of CH the cardinal $\lambda$ can be quite arbitrary,
see \cite{Je}.
We will show the next corollary which will follow immediately from 
our Main Theorem in Section 3.

\medskip

\noindent{\bf Main Corollary \ref{Theo End}}
{\em Let $A\neq 0$ be a $R$--free $R$-algebra over a countable, 
principal ideal domain  
$R$  which is not a field and let  $|A| < \lambda \leq  2^{\aleph_0}$, 
then there exists an $\aleph_1$--free $R$-module $G$  of cardinality 
$\lambda $ with $\End G = A$.
}
\medskip 

We will construct $G$ as an $A$-module and $A$ is identified with
endomorphisms acting by scalar multiplication.
If $A = R$, we derive the 
{\em existence of $\aleph_1$--free $R$-modules of cardinality 
$\aleph_1$ with $\End G = R$, a result about indecomposable 
$R$-modules} 
known only in the case $R = \Z$ from our recent paper \cite{GS1}.
The main difficulty in passing from $\Z$ to $R$ can be seen in the
local case when $R$ is a local ring with just one prime $p$, e.g. if
$R = \Z_{(p)}$ is the ring of integers localized at $p \neq 0$. 
Infinitely many primes -- by arithmetic -- provide a rigid system
(= modules with no homomorphisms $\neq 0$ between them). Hence 
homomorphisms can be restricted in their activity on $G$ by building
into $G$ a rigid system in a suitable way \cite{GS1}. Finally 
they `calm down' to scalar multiplication on $G$. 
This is no longer possible in the 
local case. The only chance we have is to utilize the existence of
sufficiently many algebraically independent elements in the $p$-adic 
completion of $\Z_{(p)}$ and this is hidden in our construction.

It may be interesting to see this result in the light of its predecessors.
The first example of an $\aleph_1$--free module which is not 
free is the Baer--Specker module $R^\omega $, which is the cartesian 
product of countably many copies of the ring $R$, known for sixty 
years; cf. Baer \cite {B} or \cite[p.94]{F1}. 
Assuming CH, this module is an example of an $R$-module of cardinality 
$\aleph_1 = 2^{\aleph_0} $. However, it is surely (by slenderness of $R$) 
a finite but not an infinite direct sum of summands $\neq 0$.
Under the same set--theoretic assumption of the
continuum hypothesis it can be shown that $A$ above can be realized 
as the endomorphism ring of an  
$\aleph_1$--free $R$-module $G$ of cardinality $\aleph_1$. The 
chronologically earlier realization theorem of this kind uses the 
weak diamond prediction principle which follows from 
$2^{\aleph_0} < 2^{\aleph_1}$. See Devlin and Shelah \cite {DS} 
for the weak diamond, Shelah \cite {S4} for the case $\End G = \Z$  and 
Dugas, G\"obel \cite {DG1} for the case $A = \End G$ and extensions to 
larger cardinals.
Using, what is called Shelah's Black Box, the existence of $\aleph_1$--free
modules $G$ with $|G| = \lambda^{\aleph_0}$  also follows from Corner, 
G\"obel \cite{CG} using Dugas, G\"obel \cite {DG2} and combinatorial 
fine tuning from Shelah \cite {S5,S6}, see also Shelah \cite[Chapter
VII]{Sh-e} and \cite{Sh-262}.   Many of the older results however do
not concentrate on the additional demand that the constructed modules
with prescribed endomorphism algebra are $\aleph_1$-free, see 
\cite{C1,G1,GM,GM1}. This of course was due to other difficulties that
had to be settled first. 
\medskip

Assuming Martin's axiom (MA) together with ZFC and  
$\aleph_2 < 2^{\aleph_0}$ any $\aleph_2$--free group $G$ of cardinality 
$< 2^{\aleph_0}$ is separable and hence has endomorphism 
ring ${\Z}$ only in the trivial case when $G =\Z$, see \cite{GS1}. 

Hence $\aleph_1$ in the Main Corollary can not be replaced by $\aleph_2$. 
This is in contrast to the result \cite {DG1} which holds in G\"odel's
universe:
All algebras $A$ as above are of the form  $A \cong \End G$  for all 
uncountable regular, not weakly compact cardinals  
$\lambda  = |G| > |A|$  such that  $G$  is  a $\lambda $--free $R$-module. 
A similar result was shown recently \cite{GS2} using the generalized
continuum hypothesis G.C.H. only.
In view of the theorem under Martin's axiom, 

{\em the existence of indecomposable $\aleph_2$--free $R$-modules of 
cardinality $\aleph_2$ and the existence
of such modules with endomorphism ring $R$, respectively, is undecidable}.

Endomorphism ring results as discussed have well-known applications using
the appropriate also well-known $R$-algebras $A$.

If  $\Gamma $  is any abelian semigroup, then we use Corner's R-algebra
$A_\Gamma $, implicitly discussed in Corner, G\"obel \cite {CG}, 
and constructed for particular $\Gamma 's$ in \cite {C2} with special 
idempotents (expressed below), with
free $R$-module structure and $|A_\Gamma | = \max \{|\Gamma |,\aleph_0\}$. If
$|\Gamma | < 2^{\aleph_0}$, we may apply the main theorem and find
a family of $\aleph_1$--free $R$-modules $G_\alpha  (\alpha  \in  \Gamma )$
of cardinality  $\aleph_1$  such that for  $\alpha ,\beta  \in  \Gamma ,$

$$
G_\alpha  \oplus G_\beta  \cong  G_{\alpha +\beta }\hbox{
and }\ \ G_\alpha  \cong  G_\beta \hbox{ if and only if }\ \ \alpha  = \beta .
$$

Observe that this induces all kinds of counterexamples to Kaplansky's test
problems for suitable  $\Gamma 's$. If we consider Corner's algebra in 
\cite {C2}, see Fuchs \cite[p.145]{F2}, then it is 
easy to see that  $A_R$ is free and  
$|A| = \aleph_0$. 
The particular idempotents in  $A$  and our main theorem provide
the existence of an $\aleph_1$--free superdecomposable $R$-module
of cardinality $\aleph_1$, which seems to be new as well. Recall that a 
group is superdecomposable if any non--trivial summand decomposes into 
a proper direct sum.

Finally, we remark that, as the reader may suspect, it is easy to replace  
$G$ in Theorem \ref{Theo End} by a rigid family of $2^\lambda $ 
such groups with only the trivial homomorphism between distinct members. 
\bigskip

\section{\sc The Construction of $\aleph_1$--free Modules}

\bigskip

{\bf a. \ The Topology}

\medskip
Let  $R$  be a countable, principal ideal domain which is not a
field, hence $R$ is reduced. 
We consider any free $R$-algebra $A$ of cardinality  $|A| < 2^{\aleph_0}$. 
In particular $A$ is torsion-free and reduced as well. 
Enumerating
$S = R\setminus \{0\} = \{ s_n : n \in \omega\}$ we obtain a descending
chain of principal ideals $q_n A$ for 
\begin{eqnarray}\label{teiler}
q_0 =1 \mbox{ and } q_{n+1} = s_n q_n^2 \hbox{ for all } n \in \omega 
\end{eqnarray}
with $\bigcap_{n \in \omega} q_n A = 0$.
The system $q_n A$ ($n \in \omega$) generates a Hausdorff topology, the 
$R$-topology on $A$. 

\bigskip
\noindent
{\bf b. \ The Geometry of a Tree and a Forest}

\medskip
Let $T = \ ^{\omega >}2$ denote the tree of all finite
branches  $v : n \rightarrow  2$, $n < \omega $, where  
$\ell (v ) = n$
denotes the length of the branch $v $. The branch of length 0  is denoted
by $\perp \  = \emptyset  \in  T$  and we also write  
$v  = (v \restr n-1) ^ \wedge v (n-1)$. 
Moreover,  $^\omega 2 = \Br(T)$  denotes all infinite
branches 
$v : \omega  \rightarrow  2$ and clearly  $v \restr n \in  T$
for all $v  \in  \Br(T)$, $n \in  \omega $. 
We often identify infinite branches $v$ with their nodes
$v = \{ v\restr n : n \in \omega \}$ which is a countable, maximal,
linearly ordered subset of $T$. Following convention we will call a
node $v\restr n$ {\it finite branch of length $n$} of the tree
$\Br(T)$. 
If $v \neq w \in \Br(T)$, then
$$br(v,w) = \inf \{ n \in \omega : v(n) \ne w(n) \}$$ 
denotes the {\it branch point} of $v$ and $w$. Hence $m = br(v,w)$ is
the largest ordinal with $v\restr m = w \restr m$.

If $C \subset\omega$, then we collect the subtree
$$T_C = \{ v \in T: \hbox{ if } e \in \ell(v) \setminus C \hbox{ then }
v(e) = 0 \}.$$
Similarly $$\Br(T_C) = \{ v \in \Br(T): \hbox{ if } e \in \omega
\setminus C \hbox{ then } v(e) = 0 \},$$
hence $v \restr n \in T_C$ for all  $v \in \Br(T_C), n \in
\omega$ and as before we omit $C$ if $C =\omega$.
Many of our arguments use a finite trunk of these trees. If $m <
\omega$, then we define
\begin{eqnarray}\label{trunk}
T^m_C = \{ \tau \in T_C: \ell(\tau) < m \}.
\end{eqnarray}

Finally let $\overline{T_C} = T_C \cup \Br(T_C)$.
\medskip

Next we use trees to build a forest.
\medskip

Let $\kappa \leq \lambda \leq 2^{\aleph_0}$ be two fixed infinite cardinals 
and     let $\kappa$ be regular and uncountable. Then we choose
a family ${\C} = \{ C_\alpha \subset\omega : \alpha < \lambda \}$ of
pairwise almost disjoint, infinite subsets of $\omega$. Let 
$T\times\alpha = \{ v\times\alpha: v \in T \}$ be a disjoint copy of the 
tree $T$ and let 
$T_\alpha = T_{C_\alpha}\times\alpha$ for $\alpha < \lambda$ be the forest
of trees (with finite branches), say 
$$T_\C = \ \stackrel{\cdot}\bigcup_{\alpha < \lambda}T_\alpha$$
and choose disjoint sets of infinite branches from $T$. We have 
$T^*_\alpha = \Br(T_{C_\alpha}) \subseteq \Br(T)$
($\alpha < \lambda$) and take a family of pairwise disjoints subsets, i.e.
$${\V} = \{ V_\alpha \subseteq T^*_\alpha: \alpha <  \lambda \} 
\mbox{ with }  |V_\alpha| = \kappa .$$
Moreover, for any $m \in \omega$ at least $\kappa$ pairs of branches
in $V_\alpha$ branch at $m$ or above.
It will be very convenient, however not necesssary to restrict to
perfect trees.
A tree is {\it perfect} if it has no isolated points (in the order
topology), i.e. every branch has an unbounded set of branch points.
It is easy to see that we may assume for the forest that
all trees $ T_\alpha \subseteq T_\C $ are perfect trees.
This additional assumption about the trees is only used in the proof
of Proposition \ref{algebra} in form of the following

\begin{observation} \label{perfect-split}
Any perfect subtree of $T$ has a subtree order isomorphic to $T$ such that
for any ordinal $n \in \omega$ there is at most one finite branch
$v\restr n$ of the subtree such that $v\restr (n+1) \neq w \restr (n+1)$ 
for some branch $w$ of this subtree.
\end{observation}

Proof: Let $T$ be a perfect tree. We will define a tree embedding 
$p :T \arr T$ such that $p(T)$ has the desired branching property.
The map $p$ is defined as the union of a chain of partial maps 
$p_i:T^i \arr T^{n_i} $. Let
$$ B(T) = \{ v\restr br(v,w) \in T: v \neq w \in T \}$$ 
be the set of all branch points in $T$. Using that $Br(T)$ has no
isolated infinite branches and K\"onig's Lemma, inductively we can 
choose a sequence of natural numbers $n_i$ ($i \in \omega$) such that 
$n_0 = 0$ and if $n_i$ is given, then $ n_{i+1}$ is the least number 
$x > n_i$ such that for any $v \in T^x$ the set
$B(v,x) = \{e : n_i < e < x,  v\restr e \in B(T) \}$ has cardinality
$|B(v,x)| > 2^{n_i}$.
Now we have enough room to extend a partial embedding
$p_i:T^i \arr T^{n_i} $ to $p_{i+1}$ in such a way 
that the branching condition of the Observation \ref{perfect-split}
holds when restricted to $p_{i+1}(T^{i+1})$.
Hence $p = \bigcup_{i \in
\omega} p_i$ is the desired tree embedding of $T$ into that perfect
tree. 
\medskip

The forest $T_\C$ of pairwise disjoint perfect trees $T_\alpha$ and the 
sequence of sets infinite branches $\V$ from $T$ which branch at `almost
disjoint sets' will form our basic geometrical objects for building modules. 
The geometry will help to distinguish elements and to carry out calculation 
in the corresponding module. In view of 
Observation \ref{perfect-split} we will assume 
\begin{eqnarray}\label{perfect}
\nonumber \mbox{ If $T_\alpha \subseteq  T_\C$ and $n \in \omega$ 
there is at most one finite branch }\\
\mbox{ $v\restr n$  such that
$v\restr (n+1) \neq w \restr (n+1)$ for some $w \in T_\alpha$.}
\end{eqnarray}

\bigskip
\noindent
{\bf c. \ The Base Module and its Completion}

\medskip
\noindent
We consider the free $A$-module 
$$B_\C = \bigoplus_{\tau \in T_\C} \tau A $$
which is a pure and dense submodule of its $R$-adic completion 
$ \BCC $ taken in the $R$-topology on $ B_\C $. The $A$-module $B_\C$ will  
be our {\em base module} and we will often omit $\C$ for convenience.
The sequence $\V$ of infinite branches is used to identify certain elements
in the completion $\BCC$. 
Any infinite branch $v \in V_\alpha$, $n<\omega$ and any $g \in \BCC$
give rise to an element $y_{vng}$ in the completion $\BCC$. Note that
\begin{eqnarray} \label{yvng}
y_{vng} = \sum_{i \geq n} \frac{q_i}{q_n} (v \restr i \times\alpha)
+ g \sum_{i \geq n} \frac{q_i}{q_n} v(i)
\end{eqnarray}
is a well--defined element of the $\A$-module $\BCC$. The reader should keep 
in mind that the {\em branch element} $y_{vng}$ connects an infinite branch
$v \in \Br(T_{C_{\alpha}}) \subseteq \Br(T)$ with finite branches from the 
disjoint tree $T_\alpha$. We will write $y_{v0g}=y_{vg}$ and often omit   
the suffix $g$ if this is clear from the context. Moreover
\begin{eqnarray}\label{n-n+1}
y_{vng} - s_n q_n y_{vn+1g} = v \restr n \times\alpha + g v(n) 
\end{eqnarray}
follows from (\ref{teiler}) and (\ref{yvng}).
We immediate obtain an equation concerning branching branches.
For $v \neq w \in V_\alpha$ and the branch point $m = br(v,w)$ we
have $v(i) = w(i) \in \{ 0,1 \}$ for $i < m$ and $v(m)-w(m) = \pm 1$.
Hence
\begin{eqnarray}
y_v = \sum_{i=0}^m q_i (v\restr i\times\alpha) 
+ g \sum_{i=0}^{m-1} q_iv(i) + q_{m+1} y_{vm+1} + g q_mv(m) 
\mbox{ and } \\
\sum_{i=0}^m q_i (v\restr i\times\alpha) + g \sum_{i=0}^{m-1} q_iv(i)=
\sum_{i=0}^m q_i (w\restr i\times\alpha) + g \sum_{i=0}^{m-1} q_iw(i) 
\mbox{ and } \nonumber
\end{eqnarray}
\begin{eqnarray} \label{tree-ele}
y_v - y_w \pm q_m g = q_{m+1} (y_{vm+1} - y_{wm+1}) \  \mbox{ for }
br(v,w) = m.
\end{eqnarray}

The special form of branch elements allows us to recognize the geometry
of the trees through the

\bigskip
\noindent
{\bf d. \ Support and Norm of Elements and Subsets of $\BCC$}

\medskip
\noindent
Elements $g$ in $\BCC$  have a natural 
{\it support} $[g] \subseteq T_\C$, the at most countable set of 
finite branches used in the sum-representation with respect to the 
$R$-adic completion of $B_\C$ defined by 
$$g = \sum_{\tau \in [g]} g_\tau , $$
where $g_\tau \neq 0$ is a unique element in $\tau \A$. 
If $\tau \in T_\C$ and $1 \in A$ we identify $\tau$ with the element
$\tau 1$ in $B$, hence $T_\C\subseteq B_\C$ and support is defined on 
$T_\C$ as well.

Let $v \in V_\alpha$ and note that in particular 
$[y_{vn0}] = [v_n\times\alpha]$, where 
$$[v_n\times\alpha] = \{ v \restr j\times\alpha: j \in \omega, 
j \geq n \}\subseteq T_\alpha.$$
This infinite part of the branch $v$ we also denote by $[v_n\times\alpha] =
[v_n]$ because it is clear that like $v$ it comes from $\alpha$.   
The notion of support trivially extends to subsets $X$ of 
$\BCC$ by taking unions $[X] = \bigcup_{x \in X}[x]$,
see also \cite{CG}.

Each element $g =\sum_{\tau \in [g]} g_\tau$ of $\B$ also has a special --
possibly empty -- subset in $[g]$, the 
$$\mbox{\em R--support \ } [g]_R = \{ \tau \in [g]: 0 \neq g_\tau \in \tau R \}.$$ 
Branch elements $g= y_{vn0}$ as in (\ref{yvng}) have R--support
$[g]_R = [g] = [v_n]$.
The support of an element immediately gives rise to its {\em norm}.

If $X \subseteq T_\C$ then 
$\|X\| = \inf\{\beta < \lambda : X
\subseteq \bigcup_{\alpha < \beta} T_\alpha \}$ denotes the {\em norm}
of $X$. If $\|X \|$ does not exist we write $\|X\| = \infty$. Moreover
$\|g\| = \|[g]\|$ denotes the norm of an element $g \in B_\C$, for
example $\|y_{vn0}\|= \alpha +1$ whenever $v \in
V_\alpha$. The following lemma is used several times.

\begin{lemma}   \label{norm}
Let $G$ be an $A$-submodule of $\BCC$ such that 
$\|G\|\leq \alpha$ and $g \in G$. If 
$$F_\alpha = \langle T_\alpha, y_{vng} : v \in V_\alpha, n < \omega
\rangle_A$$
is a submodule of $\BCC$, then $G \cap F_\alpha = Ag.$
\end{lemma}
\bigskip
Proof: Recall that $v(n) \in \{0,1\}$. Take any $n < \omega$ with 
$v(n) = 1$. From (\ref{n-n+1}) we have 
$$g = y_{vng} - s_n q_n y_{vn+1g} - (v\restr n)\times\alpha$$
and $Ag \subseteq F_\alpha \cap G$.
If $h \in F_\alpha$, then $h$ can be expressed as a linear combination 
of elements from a finite subset of $T_\alpha$ and a finite set of 
elements of the form $y_{vng}$. Using (\ref{n-n+1}) we can choose $m \in
\omega$ subject to the following conditions.
$$ h \in \langle \{g\} \cup T^m_\alpha \cup \{y_{vmg} : v \in
E \} \rangle_A \subseteq F_\alpha $$
where $E$ is a finite subset of $V_\alpha$ and
$T^m_\alpha$ as defined in (\ref{trunk}) such that 
\begin{eqnarray}
\hfill [g] \cap [v_m]  &  = & \emptyset \mbox{ for all } v \in E \\
\hfill [v_m] \cap [w_m] & = & \emptyset \mbox{ for all } v \neq w \in E 
\end{eqnarray}
We can write $h = \sum_{v \in E} a_v y_{vmg} + a_gg + t$ where 
$a_v, a_g \in A$ and $t \in \langle T^m_\alpha\rangle_A$. 
If also $h \in G$ we take any $\tau \in [v_m]$ to see that $a_v = 0$ and
similarly $t = 0$, hence $h = a_gg \in Ag$ which shows the lemma.

\bigskip
\noindent
{\bf e. \ The Desired $\aleph_1$--free Module}

\medskip
\noindent
We use these basic tools to construct the desired $R$-module.

\begin{construction} \label{constr}
Let $A$ be the free $R$-algebra over the countable principal ideal domain $R$
which is not a field  
with $|A| < \lambda \leq 2^{\aleph_0}$ and  $\kappa = |A|^++\aleph_1$.
If $\C, \V \mbox{ and } B_\C, \BCC$ are as given, then 
choose a transfinite sequence $b_\alpha$ ($\alpha < \lambda$) which
runs $\lambda$ times through the non-zero elements $b$ in $B_\C$ with
$B_\C/Ab$ $A$-free.
We define inductively $A$-submodules 
$G_\alpha \subseteq \BCC$ subject to the following 
condition for any $\alpha < \lambda$.
The sequence $G_\alpha$ is increasing,  continuous with 
$$G_0 = 0 \mbox{ and } G = \bigcup_{\alpha < \lambda}G_\alpha.$$
\begin{eqnarray} \label{gen}
G_{\alpha+1} = \langle G_\alpha \cup T_\alpha \cup \{y_{vng_\alpha}
: v \in V_\alpha, n \in \omega \} \rangle_A .
\end{eqnarray}
We also let $g_\alpha = b_\alpha$ if $b_\alpha \in G_\alpha$ with
$\|g_\alpha\| \leq \alpha $ and
$g_\alpha = 0$ otherwise.
\end{construction}

Note that $\kappa \leq \lambda \leq 2^{\aleph_0}$ is a regular cardinal.
The constructed $A$-module $G$ has visibly cardinality $\lambda$, and 
we want to show that it is $\aleph_1$-free as $A$-module.  We reserve
$G$ to denote this module for the rest of this paper.
\medskip

In view of Pontryagin's theorem we say that an $R$--module is
$\aleph_1$--free if any submodule of finite rank is contained in a free
$R$--submodule. If $R$ is a principal ideal domain, Pontryagin's theorem 
\cite[p.93, Theorem 19.1]{F1} ensures that any countably generated 
submodule is free. This gives us the following
\begin{observation} \label{Pontrj}
Let $A$ be a free $R$-algebra over a principal ideal domain $R$ and $M$ 
an $A$--module such that any finite subset is contained 
in an $A$-free and $R$-pure submodule,  then $M$ is an 
$\aleph_1$--free $R$-module, i.e. all its countably generated $R$-submodules 
are free.
\end{observation}

Next we will show that $G$ is $\aleph_1-$free. This will be the case
when $\beta = 0$ in the next proposition.

\begin{proposition} \label{free}
Let $A$ be a free $R$-algebra and $G = \bigcup_{\alpha \in \lambda} G_\alpha$ 
be the constructed $R$-module. Then $|G| = \lambda$ and $G/G_\beta$ is an
$\aleph_1-$free $R$-module for any $\beta < \lambda$.
\end{proposition}

Proof: In view of Observation \ref{Pontrj} we consider any non-empty
finite set $E \subseteq G/G_\beta$. 
Choose $\alpha < \lambda$ minimal with 
$E \subseteq G_\alpha/G_\beta$. First note that $\alpha > \beta$ must be a
successor because $E$ is a proper finite set, hence $\gamma = \alpha -1 \geq
\beta$ exists. Also note that $G_\alpha/G_\beta$ is a quotient of
$A$-modules, hence an $A$-module. By induction it is enough to show that 
\begin{eqnarray} \label{p}
E \subseteq (U+G_\beta)/G_\beta \oplus G_\gamma/G_\beta
\subseteq_* G_{\gamma +1}/G_\beta \mbox{ for some free $A$-module } 
(U+G_\beta)/G_\beta.
\end{eqnarray}
First we want to find inductively an $A$-submodule 
$U \subseteq G_\alpha$.
We note by (\ref{gen}) that $G_{\gamma +1} = G_\gamma + F_\gamma$ where
\begin{eqnarray} \label{Ggamma}
F_\gamma = \langle T_\gamma \cup  \{ y_{vng_\gamma}: v \in V_\gamma, n <
\omega \} \rangle_A
\end{eqnarray}
If $E'\subseteq G_{\gamma +1}$ is a set of representatives of the
elements in $E$, then by (\ref{Ggamma}) and (\ref{n-n+1}) 
there is a finite set 
$F \subseteq V_\gamma$, and a number $m < \omega$ such that 
$$ E' \subseteq U + G_\gamma \mbox{ \ where \ } 
U = \langle T^m_\gamma \cup \{ y_{vmg_\gamma} : v \in F \} \rangle_A$$          
Moreover we may assume that
$ [v_m] \cap [w_m] = \emptyset  \mbox{ for all } v \neq w \in F. $
A support argument shows that the defining generators of $U$ are
$A$-independent {\it modulo} $G_\gamma$, hence $U+G_\beta/G_\beta$
must be $A$-free and $G_\gamma/G_\beta \cap (U+G_\beta )/G_\beta = 0.$

Now it is easy to show that $U$ is $R$-pure in $G$ which also implies
the purity in (\ref{p}). 
If $h \in G \setminus G_{\gamma + 1}$, then an easy support argument shows
that $G_{\gamma +1}$ is pure in $G$ that is to say that $dh \not\in
G_{\gamma +1}$ for any $0 \neq d \in R$ and in particular 
$dh \not\in U$. We may suppose that $h \in G_{\gamma +1}$, 
and by the last considerations we find a finitely generated $A$-submodule
$$U' = \langle T^{m'}_\gamma \cup \{ y_{vm'g_\gamma} : v \in F' \} 
\rangle_A$$          
for some number $m'\geq m$ and finite set $F \subseteq F'\subseteq
V_\gamma$ with $h \in U'$. 
We may assume that $m'$ is chosen such that also
$$ [v_{m'}] \cap [w_{m'}] = \emptyset  \mbox{ for all } v \neq w \in F'. $$
One more support argument now shows that $U$
is a summand of $U'$, we leave it to the reader to write down a
complement of $U$ in $U'$. If $dh \in U$ for some $0 \neq d \in R$, 
then $h \in
U$ follows from $h \in U'$, which shows that $U$ is pure in $G$.

\section{\sc The Constructed Modules And Their Endomorphism Algebras}

The following Definition \ref{rigid tree} rigid tree-elements 
is the critical tool of this paper. 
The short proof of our Main Theorem \ref{Theo End},
following immediately below, is based on a Main-Lemma \ref{Main Fact} 
which indicates our strategy. Moreover, the Definition \ref{strong rigid} 
explains how to convert rigid tree-elements into algebraic content.
We think that it may help the reader if we start at the end:

The main result of this paper is the following

\begin{theorem}\label{Theo End}
If  $A$  is a free $R$-algebra over a countable, principal 
ideal domain $R$ which is not a field and 
$|A| < \lambda  \leq 2^{\aleph_0}$, then there exists an 
$\aleph_1$--free $R$-module $G$  of cardinality $\lambda$ 
with  $\End G = A.$
\end{theorem}
\medskip

Remark: $G$  will be the $A$--module constructed in (\ref{constr}) 
and we have identified $a \in  A$  with  $a\cdot id_G.$
\medskip

Proof: Let $G$ be the $A$-module from the Construction \ref{constr}.
Clearly $A \subseteq \End G$ by scalar-multiplication because $A$ acts
faithfully on $G$, and $G$ is an $\aleph_1-$free $R$-module of cardinality
$\lambda$ by Proposition \ref{free}. It remains to show that 
$\End G \subseteq A$.
\medskip

Suppose $\varphi \in \End G \setminus A$. Recall from (\ref{constr}) 
that $T_\alpha \subseteq G$ for all $\alpha < \lambda$, hence 
$B_\C \subseteq G$. Inspection of (\ref{gen}) shows that $G/B_\C$ is
torsion-free divisible. This is needed to prove that there exists
\begin{eqnarray} \label{contra}
g \in B_\C \mbox{ with } B_\C/gA \mbox{ $A$-free and } g\varphi \not\in Ag.
\end{eqnarray}
Note that $B_\C$ is a free $A$-module freely generated by some set 
$J \subseteq B_\C$. If (\ref{contra}) does not hold, then $e\varphi \in Ae$ for all
$e \in J$, say $e\varphi = a_e e$. If also $f \in J \setminus \{ e \}$ then
$f \varphi = a_f f$ and similarly $e+f$ is another basic element and
the negation of (\ref{contra}) would also give
$(e+f)\varphi =a_{e+f}(e+f)$ for some $a_{e+f} \in A$, hence
$a_e = a_{e+f} = a_f$ by independence. The element $a =a_e$ does not 
depend on $e \in J$, and $e\varphi = ae$ for all $e \in J$, hence 
$\varphi\restr B_\C = a \cdot id_{B_\C}$. The endomorphism extends 
uniquely to the $A$-module $G$ by density, and 
$\varphi = a \cdot id_G \in A$, which was 
excluded. Condition (\ref{contra}) is shown.

By Construction \ref{constr}  we can find $\alpha ' < \lambda$ such that the
element $g$ from (\ref{contra}) belongs to $G_{\alpha '}$, moreover we find
$\alpha ' < \alpha < \lambda$ with  $g = g_\alpha$, hence $G_{\alpha'}
\subseteq G_\alpha$ and $g = g_\alpha \in G_\alpha$.
In particular $g\varphi \neq 0$ by (\ref{contra}) and, since $G$ is
reduced, we find $m_0 \in \omega$ such that  
$$ g = g_\alpha, \mbox{ and } g\varphi \notin q_mG \mbox{ for all } m \geq m_0.$$
We now apply (\ref{gen})
\begin{eqnarray}
G_{\alpha+1} = \langle G_\alpha \cup T_\alpha \cup \{ y_{vn} : v \in 
V_\alpha, n \in \omega \} \rangle_A \mbox{ where } y_{vn} = y_{vng}
\end{eqnarray}
and (\ref{tree-ele}) implies for $v \neq w \in V_\alpha$ that
$$ y_v - y_w \pm q_m g = q_{m+1} (y_{vm+1} - y_{wm+1})  \mbox{ if }
br(v,w) = m. $$
We may assume $m \geq m_0$ by assumption on $\V$.
Let $t_v = y_v \varphi$  ($v \in V_\alpha \in \V$) and apply $\varphi$ to 
the last equation. We derive the
existence of a {\it family of rigid tree-elements} as defined below:

The set of elements $t_v \in G$ ($v \in T$), $T$ a subset of
$V_\alpha$ satisfying the hypothesis of (\ref{Main Fact}),
constitute a family of rigid
tree-elements for $g\varphi$, where
$$t_v - t_w \pm q_m g \varphi \in q_{m+1} G. $$
We now apply our Main-Lemma \ref{Main Fact} and obtain
$$ g \varphi \in Ag_\alpha = Ag$$
which contradicts (\ref{contra}). Hence $\varphi$ does not exist and the
Main Theorem \ref{Theo End} follows.

\medskip
We proceed with the definition of {\it rigid tree-elements}.

\begin{definition} \label{rigid tree}
Let $G$ be the $A$-module and $\kappa \leq \lambda$ be the regular
cardinal from (\ref{constr}) and $T \subseteq V_\alpha$ for some
$\alpha < \lambda$ be a set of cardinality $\kappa$. If $\beta \leq
\lambda$, a family  $\{ t_v \in G_\beta: v \in T\}$ is called a family 
of rigid tree-elements for $z \in G$ at a tree $T_\alpha$, if
$$t_v - t_w \pm q_mz \in q_{m+1}G \mbox{ for all } v \neq w \in T
\mbox{ with } br(v,w) = m.\eqno{(t)}$$
\end{definition}

\medskip

\begin{mlemma}\label{Main Fact} 
Let $G$ be the $A$-module constructed in
(\ref{constr}) and $\alpha < \lambda$, $ \beta \leq \lambda$. 
If $T$ is a set of branches from $V_\alpha$ branching above some
$m \in \omega$ with $|T| = \kappa$ and 
$\{t_v \in G_\beta: v \in T\}$ is a family
of rigid tree-elements for $ z \in G \setminus q_mG$, 
then $z \in Ag_\alpha$.
\end{mlemma}

\medskip

The proof of (\ref{Main Fact}) follows after a number of steps where we
replace $T$ above by equipotent subsets with `stronger' families of
rigid tree-elements. Our final goal is a family as in the
following definition.
\begin{definition} \label{strong rigid}
Let $m_0$ be a natural number.
We will say that the family $ \{t_v \in G_\beta : v \in T \}$ as in 
(\ref{rigid tree}) is an {\em independent} family of rigid tree-elements
for $z$ at some tree $T_\alpha$ over $m_0$ if there are a 
sequence of ordinals $\alpha_1 < \cdots < \alpha_s < \lambda$, 
$m, j^* > m_0$, a finite set $F$ with elements $a_x \in A$ for $x \in F$
and an injective map
$$\delta: T \times F \arr \ \stackrel{\cdot}\bigcup_{i\leq s}V_{\alpha_i},
\eqno{(i)}$$
such that any $t_v$ ($v \in T$) can be expressed as
$$t_v = \sum_{x \in F} a_x y_{\delta(x,v)m}. \eqno{(ii)}$$
$$\mbox{Moreover, }     \left(\bigcup_{i \leq s}\ [g_{\alpha_i}]\ \right) 
\cap [y_{\delta (x,v)m}] = \emptyset,  \eqno{(iii)}$$
$$\mbox{ the branches } \delta(x,v)\restr j^* \ (x \in F) 
\mbox{ are all distinct and independent of } v, $$
$$F = \ \stackrel{\cdot}\bigcup_{i \leq s}F_i \mbox{ and } 
\delta(T \times F_i) \subseteq V_{\alpha_i} \ (i \leq s). \eqno{(iv)}$$
\end{definition}

\medskip

For convenience we will some times omit $\alpha$ above, writing $V_i$ for
$V_{\alpha_i}$. Often we only deal with partial maps of $\delta $
say
$$ \delta_v = \delta \restr \{v \} \times F \mbox{ and write }
\delta_v: F \arr F_v = \Im (\delta_v). $$ 
In order to find an independent family of rigid tree-elements we first
concentrate on finding a weaker family which satisfies 
($w$) in (\ref{phl}) and comes from a given family of rigid tree-elements at some fixed tree.

\begin{phlemma} \label{phl}
Let $G$ be the $A$-module, $\kappa$ be the cardinal given by
(\ref{constr}), and let $\alpha < \lambda$ and $\beta \leq \lambda$.
Assume that we also have a family of elements
$$t_v \in G_\beta  \ (v \in T) \mbox{ for some subset }
T \subseteq V_\alpha \mbox{ of cardinality } \kappa $$
and elements $h_1,...,h_k \in G_\beta$ with the property that
$$ \mbox{ if } v,w \in T \mbox{ and } br(v,w) = m, \mbox{ then }
t_v - t_w \in \langle h_1,...,h_k \rangle_A +
q_mG_\beta.\eqno{(\beta)} $$
Then we find an equipotent subset $T'$ of $T$ and ordinals 
$\beta_1 < \cdots < \beta_s < \beta$ with       
$$ t_w \in \langle T_{\beta_i}, y_{vng_{\beta_i}} : v \in 
V_{\beta_i}, n < \omega \mbox{ and } i \leq s \rangle_A
\mbox{ for all } w \in T'. \eqno{(w)} $$
\end{phlemma}

\medskip
Note that ($\beta$) is a weak form of the definition of
a family of rigid tree-elements. If for the above family
$\beta = \lambda$, then by ($w$) we also find $\beta < \lambda$ such that 
$t_v \in G_\beta$ ($v \in T$) for an equipotent subfamily.

\medskip
Proof: We must collect a small `pigeon-hole' -- the right-hand side of
($w$) -- and enough `pigeons' $t_w$ to land in ($w$). There are plenty of
pigeons and we just discard all trouble makers.

The proof is by induction on $\beta \leq \lambda$.

If $\beta = 0$, then $G_0 = 0$ and it is nothing to show.

Next we assume that $\beta \leq \lambda$ is a limit ordinal. We will
distinguish three cases (a),(b) and (c) depending on the cofinality
$cf(\beta)$ of $\beta$.

\noindent
(a) \ Suppose  $cf(\beta) > \kappa.$ \\
Hence $|T| < cf(\beta)$ and note that $cf(\beta)$ is a regular cardinal.
Then we can find $\gamma < \beta$ with $h_1, ... , h_k \in G_\gamma$
and  $t_v \in G_\gamma$ for all $v \in T$. Condition ($\gamma$) in
(\ref{phl}) holds       for $\gamma < \beta$ and the induction hypothesis 
applies to finish this case.

\noindent
(b) \ Suppose  $cf(\beta) < \kappa.$\\
Note that $G_\beta = \bigcup_{i < cf(\beta)} G_{\sigma(i)}$ for some
strictly increasing, continuous sequence $\sigma(i)$ converging to
$\beta$. In this case we immediately find some $i < cf(\beta)$ and
$\gamma = \sigma(i) < \beta$ with $h_1, ... h_k \in G_\gamma$ and $T'
= \{ v \in T: t_v \in G_\gamma \}$ has cardinality $\kappa$. The
claim follows by induction like (a).

\noindent
(c) \ Suppose  $cf(\beta) = \kappa.$\\
Let     $\langle \beta_\xi \in \beta: \xi < \kappa \rangle$ be a strictly
increasing, continuous chain of ordinals converging to $\beta$, and
enumerate $T = \{ v_\xi : \xi < \kappa \}$ without repetition. We
also may assume that  $h_1, ..., h_k \in G_{\beta_0}$. Consider the
set
$$ Y = \{ \xi \in \kappa : t_{v_\xi} \in G_{\beta_\xi} \} \subseteq
\kappa.$$
Suppose first that \\
(i) \ $Y$ is a stationary subset of $\kappa$.\\
If $\xi \in Y$, then we find a smallest $G_{\beta_\delta} \ (\delta
\in \kappa)$ which contains $t_{v_\xi}$. Since $\langle \beta_\xi :
\xi \in \kappa \rangle$ and $\langle G_\nu : \nu \in \kappa \rangle$
are continuous sequences, $\delta$ must be a successor ordinal, say
$\delta = \gamma + 1$. We get a function $g: Y \arr \kappa \ (\xi \to
g(\xi) = \gamma)$ and note that $g(\xi) < \xi$ from $t_{v_\xi} \in
G_{\beta_\xi}.$ Hence $g$ is regressive and Fodor's lemma applies,
see \cite[Theorem 22, p. 59]{Je}. There is a stationary subset $X$ of
$Y$ - which must be stationary in $\kappa$ by hypothesis (i) - on
which $g$ is constant, taking some fixed value $\delta \in \kappa$.
Clearly $|X|= \kappa$ and $T'= \{t_{v_\xi} : \xi \in X \}$ and
$G_{\beta_\delta} \ (\beta_\delta < \beta)$ satisfy the induction
hypothesis ($\beta_\delta$) in (\ref{phl}). Again, the claim follows 
in case (i) by induction.\\
Finally we assume that\\
(ii) \ $Y$ is not a stationary subset of $\kappa$.\\
In this case we have to work showing that (ii) can not occur. There
is a cub $C$ in $\kappa$ with $Y \cap C = \emptyset$. Inductively we
may replace $C$ by an equipotent subset, called $C$ again, and replace
the $\beta_\xi$'s by new ones such that
$$t_{v_\xi} \in G_{\beta_{\xi+1}} \setminus G_{\beta_\xi} \mbox{ for
all } \xi \in C.$$
Note that $T' = \{ t_{v_\xi} : \xi \in C \}$ still has cardinality
$\kappa$, and apply Proposition \ref{free} to note that
$G_{\beta_{\xi+1}}/G_{\beta_\xi}$ is $\aleph_1$-free for all $\xi \in
C$~. Moreover, 
$$0 \neq t_{v_\xi} + G_{\beta_\xi} \in G_{\beta_{\xi+1}}/G_{\beta_\xi}$$ 
and there are elements $0 \neq d_\xi \in R$ such that
$$t_{v_\xi} + G_{\beta_\xi} \not\in
d_\xi(G_{\beta_{\xi+1}}/G_{\beta_\xi}). \eqno{(d)}$$
\relax From $|R| < |C|$ we can find $0 \neq d \in R$ 
with $C'= \{ \xi \in C: d_\xi = d \}$ of cardinality $\kappa$. 
Pick any $j \in \omega$ with $d \mid q_j$.

The set $X'= \{ t_{v_\xi} : \xi \in C'\}$ of cardinality
$\kappa$  must have elements $t_{v_\epsilon},
t_{v_\eta}$ for $\epsilon < \eta \in C'$ with branch point
$br(v_\epsilon, v_\eta) = m > j$. In particular $d \mid q_m$.
We derive from the hypothesis ($\beta$) that
$$t_{v_\epsilon} - t_{v_\eta} \in \langle h_1, ... , h_k \rangle_A +
q_mG_{\beta_{\eta +1}}.$$
However  $h_1, ... ,h_k \in G_{\beta_0} \subseteq G_{\beta_\eta}$ and
also $t_{v_\epsilon} \in G_{\beta_\eta}$ from $\epsilon < \eta$. Hence
$$t_{v_\eta} + G_{\beta_\eta} \in
q_m(G_{\beta_{\eta+1}}/G_{\beta_\eta}) \subseteq
d(G_{\beta_{\eta+1}}/G_{\beta_\eta})
=d_{v_\eta}(G_{\beta_{\eta+1}}/G_{\beta_\eta})$$
which contradicts ($d$) and case (ii) cannot come up.
This finishes the case of limit ordinals $\beta$.

We may assume that
$\beta = \gamma + 1$ is a successor ordinal, and the lemma holds for
$\gamma < \beta$. We also have
$$G_\beta = G_{\gamma+1} = \langle T_\gamma, G_\gamma , y_{vng_\gamma}
: v \in V_\gamma, n \in \omega \rangle_A.$$

As in (\ref{gen}) and (\ref{Ggamma}) we can write $G_{\gamma +1} =
G_\gamma + F_\gamma$ with
$$ F_\gamma = \langle T_\gamma \cup \{y_{vng_\gamma} : v \in V_\gamma,
n < \omega \}\rangle_A .$$
Obviously $\|G_\gamma\|\leq \gamma$ by construction
(\ref{constr}) of $G$ and Lemma \ref{norm} applies. We derive
$$Ag_\gamma =   G_\gamma \cap F_\gamma$$
and $G_\gamma/Ag_\gamma \oplus F_\gamma/Ag_\gamma$ is a direct sum.

If $h_1,...,h_k \in G_\beta$ and $ t_v \in G_\beta \ (v \in T)$ are
given by hypothesis, then we can write
$t_v = t_v^0 + t_v^1$ for all $v \in T$ and similarly $h_i = h_i^0 +
h_i^1 \ (i \leq k)$ with  $t_v^0, h_i^0 \in G_\gamma$ and $t_v^1,
h_i^1 \in F_\gamma$. Moreover, if $v,w \in T$ branch at $br(v,w) =
m$, then by hypothesis
$t_v - t_w \in \langle h_1,...,h_k \rangle_A + q_mG_\beta$ and 
$$(t_v^0 -t_w^0 )+ (t_v^1 - t_w^1) = \sum_{i=1}^k a_i h_i^0 + 
\sum_{i=1}^k a_ih_i^1 + q_m g'_\gamma + q_mf_\gamma$$
follows for some $g'_\gamma \in G_\gamma$ and $f_\gamma \in F_\gamma$.
Hence
$$(t_v^0-t_w^0 - \sum_{i=1}^ka_ih_i^0) - q_m g'_\gamma =
( t_w^1-t_v^1 + \sum_{i=1}^ka_ih_i^1) + q_m f_\gamma.$$
The left-hand side of the displayed equation is in $G_\gamma$ while the
right-hand side is in $F_\gamma$. The sum must be $0$ {\em modulo}
$Ag_\gamma$ by the direct sum above.
In particular
$$t_v^0 - t_w^0 \in \langle 
h_1^0,...,h_k^0, g_\gamma \rangle_A + q_mG_\gamma.$$
The induction applies for $\gamma < \beta$. We find $T' \subseteq T$ with
$|T'| = \kappa$ and $\beta_1 < ... < \beta_k < \gamma$ such that
$$t_w^0 \in \langle T_{\beta_i}, y_{vng_{\beta_i}} : v \in V_{\beta_i},
i \leq k, n \in \omega \rangle_A$$
for all $w \in T'$. Finally let $\beta_{k+1} = \gamma$ and note that
$$t_w = t^0_w+t^1_w \in 
\langle T_{\beta_i}, y_{vng_{\beta_i}} : v \in V_{\beta_i},
i \leq k+1, n \in \omega \rangle_A$$
for all $w \in T'$. This completes the induction.

\bigskip

Next we will use the Pigeon-Hole-Lemma \ref{phl} to find an independent
family of rigid tree-elements from an ordinary family of rigid
tree-elements. 
Then we are ready to prove the Main-Lemma \ref{Main Fact} which
established already the Main-Theorem \ref{Theo End}.

\begin{proposition}\label{indep}
Let $m_0$ be a natural number and $G$ be the $A$-module constructed in 
(\ref{constr}). Suppose there is a family of rigid tree-elements for 
 $0 \neq z \in G$ at some tree. Then we
can find an independent family of rigid tree-elements over $m_0$ for $z$ 
at the same tree.
\end{proposition}

Remark: The new family need not be a subfamily of the old one.

\medskip

Proof: Let $\{ t_v \in G: v \in T \}$ be the given family of rigid 
tree-elements for $z \in G$ where $T \subseteq V_\alpha$ for some 
$\alpha < \lambda$ has cardinality $\kappa$.

Hence
$$t_v - t_w \in Az + q_{m+1}G \mbox{ for all }  v \neq w \in T \mbox{
with }  br(v,w) = m$$
follows from ($t$) in Definition \ref{rigid tree}. Shrinking $T$, we
may assume $m > m_0$. Then we apply the
Pigeon-Hole-Lemma \ref{phl} for $k=1$. There is an equipotent subset
$T$ replacing $T$ and there are ordinals $\beta_1 < \cdots < \beta_s
< \lambda$ with
$$t_v \in  \sum_{i=1}^s F_i \mbox{ for all } v \in T$$
where we write [as before in (\ref{Ggamma})]
$$F_i = \langle T_{\beta_i} \cup \{ y_{wng_{\beta_i}} : w \in
V_{\beta_i}, n < \omega \} \rangle_A.$$
We replace $\beta_i$ by $i$ and put $$g_{\beta_i}=g_i, T_{\beta_i} =
T_i, V_{\beta_i} = V_i \mbox{ and } y_{wng_{\beta_i}} = y_{wn}
\mbox{ for } w \in V_i.$$
Moreover let $\BB = \{ \beta_1, ..., \beta_s \}$ and 
$\BBB = \langle T_i : i \leq s \rangle_A.$
For each $v \in T $ we now can write
$$t_v = b_v + \sum_{i=1}^s \sum_{l=1}^t
\sum_{n=0}^{m(v_{li})}a_{v_{li}n} y_{v_{li}n} \eqno{(rp)}$$
where $b_v \in \BBB , a_{v_{li}n} \in A$ and $ y_{v_{li}n}$ are
branch-elements with $v_{li} \in V_i$ depending on 
$v \in T$, $l \in \omega$ and $i \leq s$.

\medskip

Next we want to improve the representation ($rp$) by using relations
 in $G$ and discarding some of
the elements from $T$.

Note that $w\restr n \in T_i$ for $w \in V_i$ and (\ref{n-n+1}) can
be applied for $g = g_i$ to replace 
$\sum_{n=0}^{m(v_{li})}a_{v_{li}n}y_{v_{li}n}$ by multiples of $g_i$, an
element in $\BBB$ and of $y_{v_{li}m(v)}$ where $m(v)$ is larger then
the maximum of all the $m(v_{li})$'s and $m_0$, which is taken over a 
finite set of at most $st+1$ numbers.
We find new elements $a_{v_{li}m(v)}, a_{vi} \in A, b_v \in \BBB \mbox{ and }
m(v) \in \omega$ and new representations for all $v \in T$ which are
$$t_v = b_v + \sum_{i=1}^s ( \sum_{l=1}^t a_{v_{li}m(v)}
y_{v_{li}m(v)}) + a_{vi} g_i . \eqno{(newrp)}$$
Moreover we may assume, enlarging $m(v)$ for each $v \in T$ up to the
supremum of all branch points of distinct pairs $\{v_{li},v_{ji}\}$ and the
finite set $\bigcup_{i \neq j \leq s} C_i \cap C_j$ where the $C_i$'s are 
from $\C$, that
$$[y_{v_{li}m(v)}] \cap [y_{v_{nj}m(v)}] = \emptyset \mbox { for all }
(li) \ne (nj),\ 1 \leq l,\ n \leq t,\ 1 \leq i,j \leq s.$$
Also $V =\bigcup_{i \leq s}[g_i] $ is finite and
$V \cap \bigcup_{li}[y_{v_{li}m(v)}] = \emptyset$ can be obtained by
enlarging $m(v)$. This ensures the first part of ($iii$) in the Definition 
\ref{strong rigid} of an
independent family of tree-elements.
Next we apply a pigeon-hole argument to simplify {\it (newrp)} even further.
Recall that ${\rm max}\{\omega, |A|, |\BBB| \} < \kappa = |T|$.
There is a subset of cardinality 
$\kappa$ of $T$ which we denote by $T$ as well, 
with the following property.

There is a finite number of parameters  for 
$1\leq i\leq s, 1 \leq l \leq t(i)$ with elements
$$t(v,i)=t(i),\  m(v)=m > m_0, \ a_{v_{li}m(v)} = a_{li},\ a_{vi} = a_i,
\ b_v = b$$
independent of $v \in T$.
Equations {\it(newrp)} become
$$t_v = b + \sum_{i=1}^s (\sum_{l=1}^{t(i)} a_{li} y_{v_{li}m}) + a_ig_i. 
\eqno{(better)}$$
Recall that 
$v_{li} \in V_i \subseteq T_i^*,\ a_{li}, a_i \in A, b \in \BBB$. 
Let 
$$E = \{ (li) : 1 \leq i \leq s, 1 \leq l \leq t(i) \},\ 
E_v = \{v_x : x \in E \} \subseteq \ \stackrel{\cdot}\bigcup_{i\leq s}
V_i$$ 
and it is easy to verify from ({\it better}) that
$\delta_v : E \arr E_v (x \arr v_x)$
is a bijection. 

We can also choose $j_v \in \omega$ large enough and
$ > m_0$ such that the restriction map
$$(\delta_v\restr j_v): E \arr \bigcup_{i \leq s} T_i : (x \arr v_x\restr j_v)$$
is injective with image 
$(E_v\restr j_v) = \{\delta_v(x)\restr j_v : x \in E \}$.
Note that $\bigcup_{i \leq s} T_i$ is a countable set while $T$ is
uncountable. 
By a pigeon-hole argument we can shrink $T$ such that
$j_v = j^* > m_0$ and $(E_v \restr j^*) = \Delta'$ are constant for
all $v \in T$, however the finite branches in $\Delta'$ that is
$\delta_v(x)\restr j^*$ ($x \in E$) are all distinct.

\medskip

In order to show that the total map $\delta$ is injective, we
replace the old family of rigid tree-elements
$t_v$ by a new family $t_v - b - \sum_{i=1}^s a_ig_i$ ($v \in T$) 
and observe that the new family is a family of rigid tree-elements 
for $z$ as well. The new family has a better representation, we can write
$$t_v = \sum_{x \in E} a_x y_{\delta_v(x)m}. $$
The set $\{ E_v : v \in T \}$ of finite sets of infinite branches
constitutes a $\Delta$-system and the $\Delta$-Lemma applies, see 
Jech\cite[p. 225]{Je}.   
There is a new equipotent subset $T$ replacing the old $T$ such that 
$$E_w \cap E_v = \Delta  \mbox{ for all } v \neq w \in T.$$
If $F_v = E_v \setminus \Delta$ for ($v \in T$) and 
$F = E \setminus \delta_v^{-1}\Delta$, 
then the $F_v{}'$s are pairwise disjoint, hence $\delta$ is injective
on $T \times F$ and
$$d = d_v = \sum_{x \in E \setminus F} a_x y_{\delta_v(x)m} \ \ (v \in T)$$
does not depend on $v$ any more. Replacing $t_v$ again by $t_v-d$, 
we obtain a new family of rigid tree-elements for $z$ at the tree 
$T_\alpha$ with the best representation
$$t_v = \sum_{x \in F} a_x y_{\delta_v(x)m} \ \ (v \in T). \eqno{(best)}$$ 
The new family is the desired
independent family of rigid tree-elements. Recall that 
$F_v = \ \stackrel{\cdot}\bigcup_{i\leq s}(F_v \cap V_i)$ 
and the preimage of this decomposition is the
decomposition of $F$ in (\ref{strong rigid}). The proposition follows.

\bigskip
The ultimate step in proving the Main-Lemma \ref{Main Fact} is 
the following proposition. The Main-Lemma \ref{Main Fact} is now 
immediate from 

\begin{proposition}\label{algebra}
Let $G$ be the $A$-module constructed in (\ref{constr}). If there is an 
independent family of rigid tree-elements branching above some $m_0$ for 
$ z \in G \setminus q_{m_0}G$ at the tree $T_\alpha$, then $z \in A g_\alpha$.
\end{proposition}

Proof: We want to extract the arithmetical strength hidden in the given 
independent family of rigid tree-elements 
$\{t_v \in G: v \in T\}$ for $z$ at some tree $T_\alpha$ 
with $T\subseteq V_\alpha$ of cardinality $\kappa$. 
By the last Proposition \ref{indep} and Definition \ref{strong rigid} 
the elements can be expressed in the form
$$t_v = \sum_{x \in F} a_x y_{\delta_v(x)m} \mbox{ for } a_x \in A$$ 
with pairwise disjoint sets $F_v$ of 
infinite branches from $\bigcup_{i=1}^sV_i$ 
where 
$$F = \ \stackrel{\cdot}\bigcup_{i \leq s}F_i, \ 
\delta_v : F \arr F_v, \ 
\delta_v(F_i) = F_v \cap V_i \subseteq V_i \ (i \leq s),$$
$$C_i = C_{\alpha_i} \in \C, \  \
V_i = V_{\alpha_i}, \ T_i = T_{\alpha_i}, \
\ \alpha_1 < \cdots < \alpha_s < \lambda.$$ 
Let $ I = \{ \alpha_i : i \leq s \} \cup \{\alpha \}$.
Recall from (\ref{teiler}) that $\frac{q_{n+1}}{q_n} = q_n s_n$. 
Hence

\begin{eqnarray} \label{elnotin}
z \not\in \frac{q_{j}}{q_{j-1}}G = q_{j-1}s_{j-1}G \subseteq
q_{j-1} G \subseteq q_{m_0}G \mbox{ for any } j > m_0
\end{eqnarray}
by assumption on $m_0$.
Moreover, note that $F_v$ is a finite set of 
distinct branches, $\bigcup_{a \neq b \in I} C_a \cap C_b$ 
is a finite subset of $\omega$, $\bigcap_{i \in \omega}q_iG = 0$
and $\bigcap_{i \in \omega}q_iA = 0$. Also note that $F_v \cap F_w =
\emptyset$ for distinct $ v,w \in T $ by Definition \ref{strong rigid}. 
We also have an element $j^* > m_0 $ satisfying (\ref{strong rigid}).
All branches $ \delta_v(x)\restr j^*$  of length $ j^*  $
 for any $  x \in F $ are pairwise distinct but independent of $v$.
\relax From these facts it is clear that the following combinatorial conditions hold.
\begin{eqnarray} \label{distinct}
\delta_v(y) \restr j^* \neq \delta_v(x) \restr j^* 
\end{eqnarray}
\begin{eqnarray}
br(\delta_v(y),\delta_w(y)) \neq br(\delta_v(x),\delta_w(x)) \geq j^*  
\mbox{ for } v \neq w \in T \ x\neq y \in F
\end{eqnarray}
\begin{eqnarray}\label{proper}
q_{j^*-1}z \in G\setminus q_{j^*}G
\end{eqnarray}
\begin{eqnarray}\label{Ci}
 \sup\bigcup_{a\neq b \in I} C_a \cap C_b  < j^*,
\end{eqnarray}
and
\begin{eqnarray}\label{notin}
 a_x \neq 0 \Longleftrightarrow a_x     \not\in q_{j^*-1}A \mbox{ for all
} x \in F.
\end{eqnarray}
We also may assume 
\begin{eqnarray}\label{lb}
br(v,w) > j^* > m_0 \mbox{ for all } v \neq w \in T.
\end{eqnarray}
Next we will show that the branch point of any two distinct branches 
$v,w \in T$ is bounded by the branch point of some `$\delta$-pair' of branches:
\begin{eqnarray}\label{br}
br(v,w) \geq br(\delta_v(x), \delta_w(x)) \mbox{ for some } x \in F
\end{eqnarray} 
If 
\begin{eqnarray}\label{qnz}
n = br(v,w), \mbox{ then } t_v - t_w \pm q_nz \in q_{n+1}G 
\end{eqnarray}
by Definition \ref{rigid tree} ($t$) of tree-elements and if also 
$t_v - t_w \in q_{n+1}G$, then $q_n z \in q_{n+1}G$. Hence 
$z \in q_ns_nG \subseteq q_{j^*-1}s_{j^*-1}G \subseteq q_{j^*-1}G$ 
by (\ref{lb}), which contradicts (\ref{elnotin}). We have 
\begin{eqnarray}\label{divi}
t_v - t_w \in q_nG \setminus q_{n+1}G \mbox{ for } n = br(v,w).
\end{eqnarray}
On the other hand
$t_v = \sum a_x y_{\delta_v(x)m}$ and 
$t_w = \sum a_x y_{\delta_w(x)m}$ for $a_x \in A$,
hence
\begin{eqnarray}\label{dif}
t_v -t_w = \sum_{x \in F} a_x (y_{\delta_v(x)m}-y_{\delta_w(x)m}).
\end{eqnarray}
If $br(\delta_v(x),\delta_w(x)) > n$ for 
all $x \in F$, then by (\ref{tree-ele}) and the last expression 
$t_v - t_w \in q_{n+1}G$ contradicts (\ref{divi}) and (\ref{br}) follows.

\bigskip

We want to calculate (\ref{dif}) more accurate and define
$$\Br(F) = \{ br(\delta_v(x),\delta_w(x)) \in \omega : x \in F 
\mbox{ with } a_x \neq 0 \}.$$ 
Then $\mbox{min}(\Br(F)) = k = k(v,w) \leq n$ follows from (\ref{br}).

Hence
$$F'= \{ x \in F: br(\delta_v(x),\delta_w(x)) = k \}$$
is a non-empty set.

We suppose that $k < n$ for contradiction. Then $k+1 \leq n$ and $q_n G \subseteq
q_{k+1}G$ and $t_v -t_w \in q_{k+1}G$  follows from (\ref{divi}).

If $x \in F \setminus F'$, then $br(\delta_v(x),\delta_w(x)) > k$, and
$a_x(y_{\delta_v(x)m} - y_{\delta_w(x)m}) \in q_{k+1}G$ and
\begin{eqnarray}\label{equiv}
0 \equiv t_v -t_w \equiv 
\sum_{x \in F} a_x (y_{\delta_v(x)m}-y_{\delta_w(x)m}) \equiv 
\sum_{x \in F'} a_x (y_{\delta_v(x)m}-y_{\delta_w(x)m})
\mbox{ mod } \frac{q_{k+1}}{q_m}G. 
\end{eqnarray}
For any $x \in F$ we let $i(x)$ be the unique integer $i \leq s$ with
$x \in F_i$. From (\ref{distinct}) we infer
\begin{eqnarray}\label{geq} 
br(\delta_v(x), \delta_w(x)) = k \geq j^*
\end{eqnarray}
for $x \in F'$ and using  (\ref{tree-ele}) and (\ref{gen}) we can 
reduce (\ref{equiv})
further 
$$ \sum_{x \in F'}a_xq_kg_{i(x)} \in q_{k+1}G.$$
If $i(x) = i(y)$ for some $x\neq y \in F'$, then $x,y \in F'\cap F_i$
and all elements $\delta_v(x), \delta_w(x),$ $ \delta_v(y), \delta_w(y)$ 
are branches of the same tree $T_i$ and if $x \neq y$, then by 
(\ref{distinct}) the pairs of branches $(\delta_v(x), \delta_w(x))$ 
and $(\delta_v(y), \delta_w(y))$ have two distinct branch points on $T_i$
at the same level $k$. This contradicts our assumption (\ref{perfect}),  
which followed from $T_i$ being a perfect tree, hence $x = y$, 
see Observation \ref{perfect-split}. 

We have seen that $|F' \cap F_i| \leq 1$ and note that $|F' \cap F_i| = 1$
for at least one of the $i$'s, because $F' \neq \emptyset$.  
We discard all other $i$'s  and may
assume
$$F' \cap F_i = \{ x_i \} \mbox{ for all } i \leq s. $$
If $g_i = g_{i(x)} =g_{\alpha_i}$
and $a_{x_i} = a_i$, then the last displayed sum becomes
$$\sum_{i \leq s}a_iq_kg_i \in  q_{k+1}G.$$

Recall from (\ref{geq}) that $k =br(\delta_v(x),\delta_w(x)) \geq j^*$,
hence $k \not\in C_i \cap C_j$ for any $i \neq j$ and $k$ can not
be the splitting level of pairs of branches from two distinct trees 
$T_i, T_j$.
Hence splitting of branches $(\delta_v(x), \delta_w(x))$ at this level $k$
can only happen at one pair, say for the one with label $i(x) = 1$.
We can reduce the last sum expression to
$$ a_1q_kg_1 \in q_{k+1}G.$$
%for a suitable $g_i$ say $g_1$[, which in fact is $g_\alpha$].
Hence 
$$a_1g_1 \in \frac{q_{k+1}}{q_k}G = s_kq_kG \subseteq q_{j^*}G$$
by (\ref{teiler}) and (\ref{geq}). However, $g_1$ is pure in $G$,
hence $0 \neq a_1 \in q_{j^*} A$ contradicts (\ref{notin}).
It follows that  $k = n$, $F = F'$ and
\begin{eqnarray}\label{n=k}
\mbox{ all pairs } \delta_v(x), \delta_w(x) \ (x \in F) \mbox{ branch
at level } \geq n = br(v,w).
\end{eqnarray}
Note that elements $a_x (y_{\delta_v(x)m} - y_{\delta_w(x)m})$ such
that $v,w$ branches strictly above $n$ are absorbed into $q_{n+1}G$
by (\ref{tree-ele}).
As before, but now for $n = k $, it follows
\begin{eqnarray} \label{F=x}
 t_v - t_w \equiv a_x (y_{\delta_v(x)m} - y_{\delta_w(x)m})
\mbox{ mod } q_{n+1}G \mbox{ and } F = \{x\}.
\end{eqnarray}
If $\delta_v(x) = v', \delta_w(x)= w'$ and $a_x = a$ then $v',w' \in
V_{\alpha_1}$ and $v,w \in T \subseteq V_\alpha$.
The independent family of $t_v$'s for $z$ by (\ref{F=x}) simply turns
into
$$ t_v - t_w \equiv a (y_{v'm} - y_{w'm})       \mbox{ mod } q_{n+1}G.$$
The pair $v',w'$ can not branch at level $ > n$ because 
$0 \not\equiv t_v - t_w         \mbox{ mod } q_{n+1}G$ by (\ref{divi}).
Hence $v,w$ and $v',w'$ branch at the same level $n >j^*$ by
(\ref{n=k}). Either $\{v,w\} = \{v',w'\}$ or the pairs are different.
In the second case branching of two distinct pairs of branches at such a 
high level $n$ can only happen at the `same tree'. 
Hence in either case we must have
$$\alpha_1 = \alpha \mbox{ and also } g_1 = g_\alpha$$ 
by (\ref{notin}). Using (\ref{tree-ele}), (\ref{gen}) and $(t)$ we have
$$0 \not\equiv q_nz \equiv aq_ng_\alpha \mbox{ mod } q_{n+1}G.$$
As before, we derive  
$z \equiv a g_\alpha \mbox{ mod } \frac{q_{n+1}}{q_n} G = q_n s_nG$ from 
(\ref{teiler}).
The set $T$ of infinite branches has size $\kappa$ and hence its
branches  split at arbitrarily large level. Choose any sequence of 
pairs ($v,w$) of branches from $T$ with branch points 
converging to infinity and note that $G$ is $\aleph_1$--free by
Proposition \ref{free}, hence $G$ is reduced. We derive  that
$z - a g_\alpha \in \bigcap_{n \in \omega} q_ns_nG = 0,$ hence 
$z = a g_\alpha \in A g_\alpha$ as required.

\bigskip
Proof of the Main-Lemma \ref{Main Fact}: The family of rigid
tree-elements given by (\ref{Main Fact}) can be traded into an
independent family of rigid tree-elements over the given number $m$
and the same tree by Proposition \ref{indep}. Now the assumptions
for Proposition \ref{algebra} are satisfied for $m=m_0$ and that $z$.
Hence the conclusion of (\ref{Main Fact}) follows from
(\ref{algebra}).

\noindent
R\"udiger G\"obel \\
Fachbereich 6, Mathematik und Informatik \\
Universit\"at Essen, 45117 Essen, Germany \\
{\small e--mail: R.Goebel@Uni-Essen.De}\\
and \\ 
Saharon Shelah \\
Department of Mathematics\\ 
Hebrew University, Jerusalem, Israel \\
and Rutgers University, Newbrunswick, NJ, U.S.A \\
{\small e-mail: Shelah@math.huji.ae.il}


\medskip

\begin{thebibliography}{}

\bibitem{B} 
 {\sc R. Baer:} {\em Abelian groups without elements of finite order}, 
 Duke Math. J. {\bf 3} (1937), 68--122.

\bibitem{C1} 
{\sc A.L.S. Corner:} {\em Every countable reduced torsion--free ring is an 
endomorphism ring}, Proc. London Math. Soc. (3) {\bf 13} (1963), 687--710.

\bibitem{C2}
{\sc A.L.S. Corner:} {\em Additive categories and a theorem of W.G. Leavitt},
Bull. Amer. Math. Soc. {\bf 75} (1969), 78--82.

\bibitem{CG}
{\sc A.L.S. Corner} and {\sc R. G\"obel:} {\em Prescribing endomorphism 
algebras, a unified treatment}, Proc. London Math. Soc. (3) {\bf 50} (1985), 
447--479.

\bibitem{DS}
{\sc K. Devlin} and {\sc S. Shelah:} {\em A weak version of $\diamondsuit$ 
which follows from $2^{\aleph_0} < 2^{\aleph_1}$}, Israel J. Math. {\bf 29} 
(1978), 239--247.

\bibitem{D}
{\sc M. Dugas:} {\em Fast freie abelsche Gruppen mit Endomorphismenring $\Z$},
Journal of Algebra {\bf 71} (1981), 314-321.


\bibitem{DG1}
{\sc M. Dugas} and {\sc R. G\"obel:} {\em Every cotorsion--free ring is an 
endomorphism ring}, Proc. London Math. Soc (3) {\bf 45} (1982), 319--336.

\bibitem{DG2}
{\sc M. Dugas} and {\sc R. G\"obel:} {\em  Every cotorsion--free algebra is 
an endomorphism algebra}, Math. Zeitschr. {\bf 181} (1982), 451--470.

\bibitem{DG3}
{\sc M. Dugas} and {\sc R. G\"obel:} {\em  On radicals and products}, 
Pacific J. Math. {\bf 18} (1985), 70--104.

\bibitem{Ed}
{\sc K. Eda:} {\em Cardinal restrictions for preradicals}, pp.277--283 in 
{\em "Abelian Group Theory"}, Contemporary Math. {\bf 87}, Providence 1989.

\bibitem{E1}
{\sc P.C. Eklof:} {\em On the existence of $\kappa $--free abelian groups}, 
Proc. Amer. Math. Soc. {\bf 47} (1975), 65--72.

\bibitem{E2}
{\sc P.C. Eklof:} {\em Set theoretic methods in homological algebra and 
abelian groups}, Les Presses de l'Universit{\'e} de Montr{\'e}al, 
Montreal 1980.


\bibitem{EM1}
{\sc P.C. Eklof} and {\sc A.H. Mekler:} {\em On 
constructing indecomposable groups in $L$}, Journal of Algebra {\bf 49} (1977), 
96-103. 


\bibitem{EM}
{\sc P.C. Eklof} and {\sc A.H. Mekler:} {\em Almost free modules, 
set--theoretic methods}, North--Holland, Amsterdam 1990.

\bibitem{ES}
{\sc P.C. Eklof} and {\sc S. Shelah:} {\em On Whitehead modules}, Journal 
of Algebra {\bf 142} (1991), 492-510.


\bibitem{F1}
{\sc L. Fuchs:} {\em Infinite abelian groups, Vol. I}, Academic Press, New 
York, 1970.

\bibitem{F2}
{\sc L. Fuchs:} {\em Infinite abelian groups, Vol.II}, Academic Press, New 
York, 1973.

\bibitem{G1}
{\sc R. G\"obel:} {\em An easy topological construction for realizing 
endomorphism rings}, Proc. R. Irish Acad. {\bf 92 A} (1992) 281-284. 


\bibitem{GM}
{\sc R. G\"obel} and {\sc W. May:}  {\em Independence in completions 
and endomorphism algebras}, Forum Mathematicum {\bf 1} (1989), 215-226.

\bibitem{GM1}
{\sc R.G\"obel} and {\sc W. May:}  {\em Four submodules suffice 
for realizing algebras over commutative rings}, J. Pure Appl. Algebra 
{\bf 65} (1990), 29-43.


\bibitem{GS1}
{\sc R. G\"obel} and {\sc S.Shelah:} On the existence of rigid 
$\aleph_1$-free abelian groups of cardinality $\aleph_1$, pp. 227--237 in
{\em Abelian Groups and Modules} 
Proceedings of the Padova Conference 1994, 
Mathematics and Its Applications, Vol. {\bf 343}, 
Kluwer Academic Publ., London  1995


\bibitem{GS2}
{\sc R. G\"obel} and {\sc S.Shelah:} G.C.H. implies the existence of 
many rigid almost free abelian groups, pp. 253--271 in
{\em Abelian Groups and Modules} 
Proceedings of the international conference at Colorado Springs 1995, 
Lecture Notres in pure and appl. mathematics, Vol {\bf 182}, Marcel Dekker,
New York 1996



\bibitem{Gr}
{\sc Phillip Griffith:} {\em $\aleph_n$--free abelian groups}, 
Quart. J. Math. (2) {\bf 23} (72), 417--425.

\bibitem{Hg1}
{\sc G. Higman:} {\em Almost free groups}, Proc. London Math. Soc. {\bf 1} 
(1951), 184--190.

\bibitem{Hg2}
{\sc G. Higman:} {\em Some countably free groups}, pp. 129--150 in 
Proceedings "Group Theory", W. de Gruyter Publ. Singapore 1991.

\bibitem{Hi}
{\sc P. Hill:} {\em New criteria for freeness in abelian groups II}, Trans. 
Amer. Math. Soc. {\bf 196} (1974), 191--201.

\bibitem{Je}
{\sc T. Jech:} {\em Set Theory}, Academic Press, New York 1978.

\bibitem{MS}
{\sc M. Magidor} and {\sc S. Shelah:} {\em When does almost free imply free?
 (for groups, transversals, etc.)}, Journ. Amer. Math. Soc. {\bf 7} (4)
 (1994), 769--830.


\bibitem{Me}
{\sc A.H. Mekler:} {\em How to construct almost free groups}, Can. J. Math. 
{\bf 32} (1980), 1206--1228.

\bibitem{R}
{\sc J. Rotman:} {\em Homological Algebra}, Academic Press, New York 1979.


\bibitem{S0}
{\sc S. Shelah:} {\em A compactness theorem for singular cardinals, 
free algebras, Whitehead problem and transversals}, Israel Journal of Math., 
{\bf 21} (1975), 319-349.


\bibitem{S1}
{\sc S. Shelah:} {\em On successors of singular cardinals}, Logic 
Colloquium '78 {\bf 97} (1978) of Stud. Logic Foundations Math., 357-380, 
North Holland, Amsterdam-New York.


\bibitem{S2}
{\sc S. Shelah:} {\em On uncountable abelian groups}, Israel J. Math. 
{\bf 32} (1979), 311--330.



\bibitem{S3}
{\sc S. Shelah:} {\em Whitehead groups may not be free, even assuming
CH, II,} Israel J. Math. {\bf 35} (1980), 257 -- 285.


\bibitem{S4}
{\sc S. Shelah:} {\em On endo--rigid strongly $\aleph_1$--free abelian groups 
in $\aleph_1$}, Israel J. Math. {\bf 40} (1981), 291--295.



\bibitem{S5}
{\sc S. Shelah:} {\em A combinatorial principle and endomorphism rings. I: On
$p$--groups}, Israel J. Math. {\bf 49} (1984), 239--257.



\bibitem{S6}
{\sc S. Shelah:} {\em A combinatorial theorem and endomorphism rings of
abelian groups, II}, pp. 37 - 86 in {\em Abelian groups and modules} 
CISM Lecture Notes, Vol. {\bf 287}, Springer Wien-New York 1984.

\bibitem{S7}
{\sc S. Shelah:} {\em Incompactness in regular cardinals}, Notre Dame 
Journal of Formal Logic, {\bf 26} (1985), 195-228.



\bibitem{S8}
{\sc S. Shelah:} {\em A note on $\kappa$-freeness of abelian groups},  
pp. 260 - 268 in {\it Around classification theory of models}, Lecture 
Notes in Math. Vol. {\bf 1182}, Springer, Berlin, 1986.

\bibitem{Sh-262}
{\sc S. Shelah:} {\em The number of pairwise non-elementary embeddable
models,} Journal of Symbol. Logic, {\bf 54} (1989), 1431 - 1455. 
\bibitem{Sh-e}
{\sc S. Shelah:} {\em Non Structure Theory,} Oxford University Press,
1996.


\bibitem{SS} 
{\sc S. Shelah} and {\sc Z. Spasojevic:} {\em A forcing axiom on strongly 
inaccessible 
cardinals making all uniformizations for all ladder systems lying 
on a fixed stationary set and applications for abelian groups,} manuscript
No. 587.


\end{thebibliography}
\end{document}